\documentclass[10pt,journal]{IEEEtranTCOM}

\usepackage[english]{babel}
\usepackage[usenames]{color}
\usepackage[cp1250]{inputenc}
\usepackage{amsfonts}
\usepackage{amsthm}
\usepackage{graphicx}
\usepackage{epsfig}
\usepackage{mathrsfs}
\usepackage{amsmath}
\usepackage{algorithm}
\usepackage{algorithmic}

\theoremstyle{plain}

\begin{document}
\newcommand{\bea}{\begin{eqnarray}}
\newcommand{\eea}{\end{eqnarray}}
\newcommand{\be}{\begin{equation}}
\newcommand{\ee}{\end{equation}}
\newcommand{\beas}{\begin{eqnarray*}}
\newcommand{\eeas}{\end{eqnarray*}}
\newcommand{\bs}{\backslash}
\newcommand{\bc}{\begin{center}}
\newcommand{\ec}{\end{center}}
\def\SC {\mathscr{C}}

\title{Improving Pyramid Vector Quantizer\\with power projection}
\author{\IEEEauthorblockN{Jarek Duda}\\
\IEEEauthorblockA{Jagiellonian University,
Golebia 24, 31-007 Krakow, Poland,
Email: \emph{dudajar@gmail.com}}}
\maketitle

\begin{abstract}
Pyramid Vector Quantizer (PVQ) is a promising technique especially for multimedia data compression, already used in Opus audio codec and considered for AV1 video codec. It quantizes vectors from Euclidean unit sphere by first projecting them to $L^1$ norm unit sphere, then quantizing and encoding there. This paper shows that the used standard radial projection is suboptimal and proposes to tune its deformations by using parameterized power projection: $x\to x^p/\|x^p\|$ instead, where the optimized power $p$ is applied coordinate-wise, getting usually $\geq 0.5\,  dB$ improvement comparing to radial projection.
\end{abstract}
\section{Introduction}
Vector quantization is seen as a promising direction for improving data compression, especially for multimedia data. Perceptual Vector Quantization is already used in Opus audio codec~\cite{opus} and is considered for AV1 video codec~\cite{daala} due to conserving the spectral envelope of audio signal and energy in video. It normalizes the input and prediction vectors, then performs Hausholder reflection to fix position of the input vector, then encodes the angle from prediction, and finally quantizes, normalizes and encodes the resulting difference as a vector from unit Euclidean sphere.

\begin{figure}[t!]
    \centering
        \includegraphics{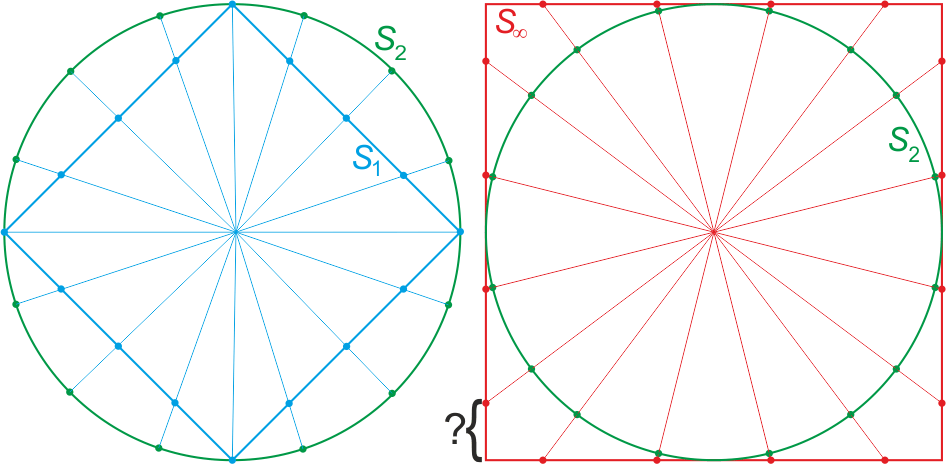}
        \caption{While in practical applications the most important seems quantization of vectors from Euclidean sphere $S_2$, direct quantization of spheres $S_1$ and $S_\infty$ is much simpler. Hence, a natural approach to perform quantization of $S_2$ is first going to $S_1$ or $S_\infty$ using some bijection (radial projection here), perform quantization, and then go back to $S_2$ using inverse bijection. A standard choice for this bijection is just radial projection $x\to x/\|x\|_s$. However, it leads to a relatively nonuniform lattice on $S_2$, which will be improved here by tuning parametrized family of power projections: $x\to x^p/\|x^p\|_s$. As in PVQ, we will focus on $S_1$ sphere approach here, which is more costly to address (enumerative coding) than $S_\infty$ (direct coding), but leads to essentially lower MSE (usually $\geq 10 \%$) as pyramidal lattice of $S_1$ is intuitively closer to sphere packing than $\mathbb{Z}^{L-1}$ lattice for $S_\infty$. However, the $S_\infty$ approach might be useful for low number of bits where PVQ has poor performance. The question mark represents its additional freedom.}
       \label{spheres}
\end{figure}

It uses Fischer's 1986 Pyramid Vector Quantizer (PVQ)~\cite{PVQ} to encode this vector from unit Euclidean sphere. It first applies radial projection to $L^1$ norm unit sphere, like presented in Fig. \ref{spheres}, then quantizes it with a uniform pyramidal lattice and calculates the number of this point in the space of all possibilities (enumerative coding). This paper argues that using radial projection is suboptimal and proposes to tune its deformation by optimizing parameter of the introduced power projection - getting inexpensive but essential improvement, which reduces mean-square error (MSE) usually by more than $10\%$.

\section{Improving Pyramid Vector Quantizer}

PVQ encoder starts with $L$ dimensional vector $x\in S_2$ from Euclidean unit sphere:
$$ S_s:=\{x\in \mathbb{R}^L: \|x\|_s=1\},\quad S_s^+:=S_s\cap (\mathbb{R}^+\cup\{0\})^L, $$
\be\textrm{where}\qquad \|x\|_s=\left(\sum_i |x_i|^s\right)^{1/s}.\ee
\noindent It first performs radial projection into $L^1$ norm unit sphere $S_1$:
\be y=P_1(x)\qquad\textrm{where}\qquad P_s(x):=x/ \|x\|_s. \ee
\noindent Then it approximates $y$ as the closest point $\tilde y=q_{LK}(y)$ from
\be S(L,K):=\left\{\left(\frac{\tilde{y_i}}{K}\right)_{i=1..L}:\tilde{y_i}\in\mathbb{Z},\ \sum_i |\tilde{y_i}|=K\right\}\subset S_1 \ee
\noindent  for a chosen parameter $K\in\mathbb{N}$ determining precision and bit cost. Instead of formally defining this subtle "closest point" condition for $q_{LK}$, there will be later presented implementation used for calculating it in the benchamarks.

\begin{figure}[t!]
    \centering
        \includegraphics{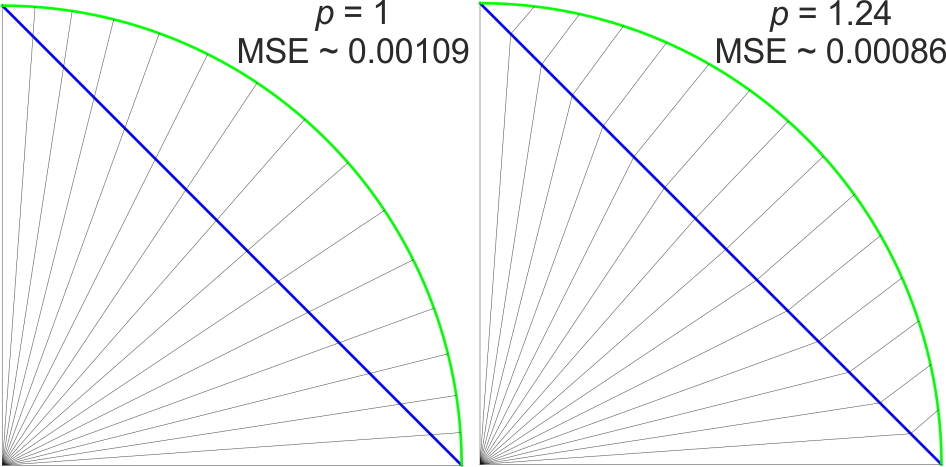}
        \caption{Two-dimensional case ($L=2,\ K=15$) for standard radial projection (left) and power projection (right) using $p=1.24$ chosen to minimize MSE - by more uniformly distributing the corresponding points. These points can be evenly distributed by just choosing trigonometric formula (\ref{even}) here, however, it does not directly generalize to higher dimensions, where the optimization criterium is also much more complex. Another suggestion from the above figure is just to use orthogonal instead of radial projection, however, again it does not generalize to higher dimensions.}
       \label{l2k15}
\end{figure}
Finally, there is used recurrent combinatorial formula to assign $\tilde y$ a number in enumeration of all points from $S(L,K)$. This final number is stored in the compressed file (fractional bits should be handled here). Decoder first decodes $\tilde y\in S(L,K)$, then performs radial projection back to Euclidean sphere $S_2$, getting: $\tilde x = P_2(\tilde y)$.

\subsection{Optimal projection}
While the above choice of radial projection might seem intuitively natural, it is not necessarily the optimal one. The standard optimization criterium is minimization of mean-square error:
\be \textrm{MSE}= \frac{1}{\textrm{volume}(S_2)} \int_{S_2}\|x-\tilde x\|_2^2\ dx,\ee
\noindent which intuitively should be minimized by trying to maintain uniform distribution of $S(L,K)\subset S_1$ while projecting it to $S_2$.
Let us first try to understand the continuous limit $(K\to\infty)$ optimization criterium for such function $f:S_1\to S_2$, allowing for differential formulation. Its Jacobian matrix (symmetric):
$$D:=[f_{i,j}]_{ij}\qquad \textrm{where}\qquad f_{i,j}:=\partial f_i/\partial x_j$$
\noindent determines local linear behavior. It performs independent scalings in the directions of eigenvectors by corresponding (real nononegative) eigenvalues. Mean-square distance scales with sum of squares of these scaling factors (eigenvalues). This sum does not depend on the chosen base: is square of the Frobenius norm: $Tr(D D^T)$. Finally, multiplying by the unit volume of our transformation, which is determinant of Jacobian matrix (assume it is positive here), the ideal transformation $f$ would be the one minimizing:

\be \int_{S_1} \det\left([f_{i,j}]_{ij}  \right)\ \sum_{ij} (f_{i,j})^2 \ dx_1\ldots dx_L \ee
\noindent Using the calculus of variations, we can transform this minimization problem into solving of Euler-Lagrange PDE for $f$:
\be \forall_i\ \frac{d}{dx_i} \frac{\partial \left(\det\left([f_{i,j}]_{ij}  \right)\ \sum_{ij} (f_{i,j})^2 \right)}{\partial f_{i,j}} =0. \ee

Unfortunately, finding such optimal $f:S_1\to S_2$ seems an extremely difficult task (still operating on minors). Even if analytical formula exists an will be found, it might turn out too costly to use in data compressors. More importantly, it assumes $K\to \infty$, making it not necessarily optimal for a finite $K$. Hence, there is a need for a practical approximation, like some inexpensively invertible parametric family of functions - which parameters will be optimized empirically.

\begin{figure}[t!]
    \centering
        \includegraphics{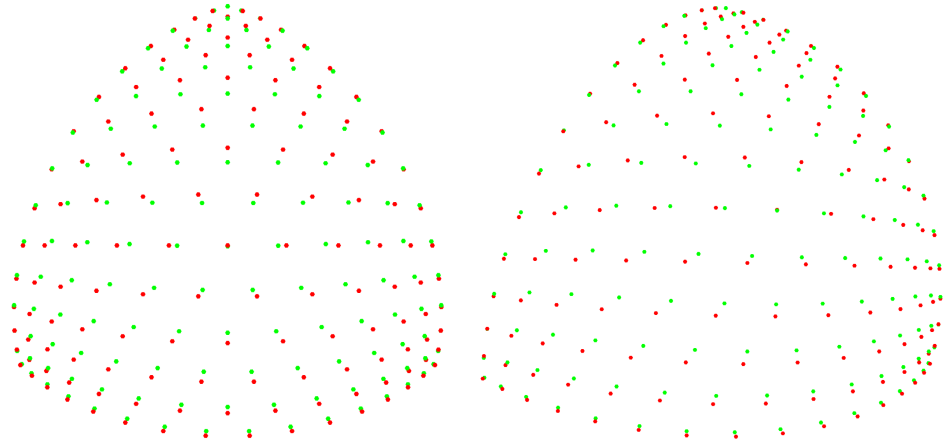}
        \caption{Two viewpoint perspectives on projecting points from $S_1^+$ for $L=3,\ K=15$ into $S_2^+$. The red points used standard radial projection, the green ones were chosen to minimize $\sum_e \|x^{e_1}-x^{e_2}\|_2^4$ over all pairs of neighboring vertices, which is some approximate of minimizing MSE for quantization using these points. As in Fig. \ref{l2k15}, the main intuition is that the priority is to condensate the central points (reduce distances), at cost of dilution near the boundaries. Power projection was arbitrarily chosen as inexpensively invertible parameterized family of functions allowing to perform similar condensation-dilution.}
       \label{2d}
\end{figure}

To gain intuitions what we really expect from such function, let us start with looking at low dimensions $L=2,3$, presented in examples in Fig. \ref{l2k15} and \ref{2d}. For $L=2$ we can easily find analytical formula uniformly distributing these points:
\be S^+_1 \ni (x_1,x_2)\to (\sin(x_1\ \pi/2),\sin(x_2\ \pi/2))\in S^+_2 \label{even} \ee
However, for $L=3$ the situation is much more complex. The green points in Fig. \ref{2d} represent solution of minimization of $\sum_e \|x^{e_1}-x^{e_2}\|_2^4$ over all pairs of neighboring vertices, which is some approximation of MSE (square distance times square for area). We see that comparing to red points, which were obtained by just radial projection, the optimized green ones are more condensed in the center (lower errors), at cost of dilution near the boundaries.

\subsection{Power projection and benchmarks}

Such function (projection) should be chosen such that both it and its inverse have very low computational cost, for example basing on some coordinate-wise operation which can be put into a table. A natural candidate is the following function, referred here as \emph{power projection}:

\be P^p_s(x):=P_s((|x_i|^p\, \textrm{sgn}(x_i))_{i=1..L})\ee

The PVQ process is modified by replacing the original radial projection with this power projection:
$$ y=P_1^{1/p}(x)\quad\to\quad \tilde y=q_{LK}(y)\quad\to\quad \tilde x=P_2^p(\tilde y) $$
\noindent Without quantization: for $\tilde y=y$, we would have $x=\tilde x$. For $p=1$ we get standard PVQ, for $p>1$ we get as required: condensation in the center of $S_2^+$, dilution at its boundaries. For $p<1$ we would get the opposite behavior.

The choice of the optimal parameter $p$ is a difficult question, it varies with $L$ and $K$. However, the optimal one is usually close to $1.3$, and its small changes in $[1.2, 1.4]$ range have usually nearly negligible impact on the MSE.

\subsection{Benchmarks}
\begin{figure}[t!]
    \centering
        \includegraphics{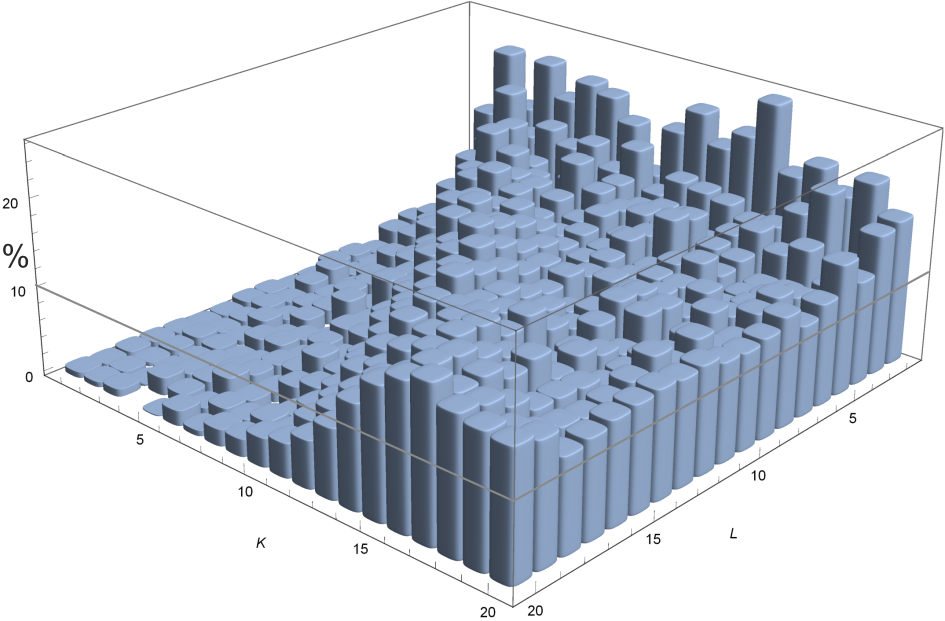}
        \caption{Numerical results for percentage reduction of MSE for power projection comparing to standard radial projection for $L=2,\ldots,20$ and $K=1,\ldots,20$. The largest found difference was $26\% (1.3\, dB)$ for $L=2,\ K=15$, presented in Fig. \ref{l2k15}. For higher dimensions $L$ the general behavior seems to be the following. For $K<L/2$ the reduction is nearly negligible (flat region on the left), then starts growing with $K$ up to $16-17\%$ ($\approx 0.8\, dB$) for $K$ slightly smaller than $L$. For larger $K$ it drops and seems to stabilize at $10-13\%$ ($\approx 0.5\, dB$). The power $p$ was optimized individually for every case, but generally it was close to $1.3$. Some further improvements could be probably obtained by choosing a function with more parameters than power projection.}
       \label{comp}
\end{figure}

The test results are presented in Fig. \ref{comp}. They were obtained by calculating mean of
\be \textrm{MSE}_p:=\textrm{average of}\ \|x-P^p_2(q_{LK}(P_1^{1/p}(x)))  \|_2^2 \ee
over 10000 random initial points $x\in S_2$, independently for every $(L,K)$ and for all $p=1$ to $1.5$ with step 0.01. The $p$ with the lowest MSE was finally chosen.  Then the shown percentage improvement is $100(1-\textrm{MSE}_p /\textrm{MSE}_1)$. It can be translated into gain in decibels as $10\log_{10}(\textrm{MSE}_1/\textrm{MSE}_p)$.

The exact choice of $q_{LK}$ is subtle - the following Mathematica function was used for quantizing the absolute values of coordinates (\verb"va"$\in S_1^+$):

\begin{footnotesize}
\begin{verbatim}
  quant[va_, k_] := (
  vk = k*va; vr = Round[vk]; kr = Total[vr];
  If[kr != k,     (* repair quantization *)
   If[k > kr, (* sort by differences: *)
     dif = vr - vk; ord = Ordering[dif];
     Do[vr[[ord[[i]]]]++, {i, k - kr}],
     (* Sign[] prevents reducing 0: *)
     dif = vk - vr - Sign[vr]; ord = Ordering[dif];
     Do[vr[[ord[[i]]]]--, {i, kr - k}]
     ]]; vr/k )
\end{verbatim}
\end{footnotesize}

\subsection{The case of low $K$}
As discussed in Fig. \ref{comp}, the real income from tuning starts with $K>L/2$. For lower $K$ the concept behind PVQ: of using a lattice on $S_1$ becomes degenerated - while only a small fraction of coordinates can have nonzero values, for each two coordinates there are used multiple points $(i,K-i)/K$. This imbalance suggests suboptimality.

Therefore, PVQ was compared with some trivial vector quantizers, like just remembering signs of all coordinates, what costs $L$ bits. For example for $L=15$ it leads to MSE $\approx 0.24$. In comparison, we need $\approx 15.06$ bits for $L=15$, $K=4$ PVQ, which gives much larger MSE $\approx 0.47$, what means $\approx 2.9\, dB$ loss. Another trivial vector quantizer: remember signs and position of the highest absolute value coordinate, what costs $L+\lg(L)$ bits, is slightly better than $K=6$ here.

The lesson here is that PVQ degenerates for $K<L/2$ cases - performance can be often significantly improved by replacing with a different quantizer. Unfortunately, optimizations here may cost the parameter flexibility of PVQ - probably requires individual treatment for different targeted numbers of bits. However, some flexibility can be gained by using $S_\infty$ sphere instead of $S_1$, like suggested in Fig. \ref{spheres}. Storing only signs is one of its degenerated cases, and generally (beside additional deformation) it has more freedom of choosing the lattice, like if the corners should be used, otherwise what is the minimal distance to corner.

\section{Conclusions and further perspectives}
The paper discussed suboptimality of using radial projection in standard PVQ, and suggested to use still suboptimal but essentially better: power projection. Depending on its tuned parameter, it usually leads to $0.5\, dB$ or larger gain, which among others can improve compression ratio of Opus audio codec (would require modifying the standard), and can help with adaptation and performance for AV1 video codec.

There have remained many open questions, like details of choosing the parameter $p$ and accuracy of tables for calculating power in implementation. Another question is finding a better function than the arbitrarily chosen power projection, for example parametrized by two parameters instead of a single one. A natural approach here is combining radial projection with some bijection inside the $S_1^+$ simplex, for example by performing some nonlinear 1D transformation inside the line segment between the center and boundary of this simplex. Such transformation can be additionally parameterized by some angle of this line segment. There is generally a large space of possibilities to test here.

Very interesting but also difficult is improving the theoretical understanding of this problem, especially the theoretical limitations for such applied deformations: how many percents can MSE be further reduced?

There was also discussed very weak performance of PVQ for low $K$, especially for $K<L/2$. It can be essentially improved by choosing a different quantizer, but optimizing the details require further work. Using $S_\infty$ sphere instead of PVQ $S_1$ sphere might be a promising direction here, still maintaining flexibility.

Finally, an essential issue of PVQ as a candidate for AV1 is extremely high computational cost. A possible solution is replacing costly combinatorial enumerative coding with a fast accurate entropy coder, which allows to provide a similar compression ratio using approximated probability distribution of symbols. It would also naturally handle the issue of storing fractional bits for optimality of PVQ.

\bibliographystyle{IEEEtran}
\bibliography{cites}
\end{document}